\documentclass{amsart}
\usepackage{amsmath,amssymb}
\newtheorem{theorem}{Theorem}[section]
\newtheorem{lemma}[theorem]{Lemma}
\newtheorem{proposition}[theorem]{Proposition}
\newtheorem{corollary}[theorem]{Corollary}

\newtheorem{remark}[theorem]{Remark}
\newtheorem{example}[theorem]{Example}

\begin{document}
\title{Pure-injective hulls of modules over valuation rings}
\author{Fran\c cois Couchot}
\address{Laboratoire de Math\'ematiques Nicolas Oresme, CNRS UMR
  6139,
D\'epartement de math\'ematiques et m\'ecanique,
14032 Caen cedex, France}
\email{couchot@math.unicaen.fr} 

\begin{abstract}If $\widehat{R}$ is the pure-injective hull of a valuation
  ring $R$, it is proved that $\widehat{R}\otimes_RM$ is the
  pure-injective hull of $M$, for every finitely generated $R$-module
  $M$. Moreover $\widehat{R}\otimes_RM\cong \oplus_{1\leq k\leq
  n}\widehat{R}/A_k\widehat{R}$, where $(A_k)_{1\leq k\leq n}$ is the
  annihilator sequence of $M$. The pure-injective hulls of uniserial
  or polyserial modules are also investigated.
Any two pure-composition series of a countably generated polyserial module are isomorphic.
\end{abstract}
\maketitle
\bigskip
The aim of this paper is to study pure-injective hulls of modules over
valuation rings. If $R$ is a valuation domain and $S$ a maximal
immediate extension of $R$, then, in \cite{War69}, Warfield proved that
$S$ is a pure-injective hull of $R$. Moreover, for each finitely
generated $R$-module $M$, he showed that $S\otimes_RM$ is a
pure-injective hull of $M$ and a direct sum of $\mathrm{gen}\ M$ indecomposable
pure-injective modules. We extend this last result to every valuation
ring $R$ by replacing $S$ with the pure-injective hull $\widehat{R}$ of
$R$. As in the domain case $\widehat{R}$ is a faithfully flat module.
Moreover, for each $x\in\widehat{R}$ there exist $r\in R$ and $y\in 1+P\widehat{R}$ such that $x=ry$.
This property allows us to prove most of the main results of this paper.
We  extend results obtained by Fuchs and Salce on 
pure-injective hulls of uniserial modules over valuation domains
(\cite[chapter XIII, section 5]{FuSa01}). We show that the length of any
pure-composition series of a polyserial module $M$ is its Malcev rank
$\mathrm{Mr}\ M$ and its pure-injective hull $\widehat{M}$ is
  a direct sum of $p$ indecomposable pure-injective modules,
  where $p\leq\mathrm{Mr}\ M$. But it is possible to have
  $p<\mathrm{Mr}\ M$ and we prove that the equality  holds for all $M$
  if
  and only if $R$ is maximal (Theorem~\ref{T:malcev}). This result is a consequence of 
the fact  that $R$ is maximal if and only if $R/N$ and $R_N$ are maximal, where $N$ is 
the nilradical of $R$ (Theorem~\ref{T:maximal}). If $U_1,\dots,U_n$ are the factors of a pure-composition
series of a polyserial module $M$ then the collection $(\widehat{R}\otimes_RU_k)_{1\leq k\leq n}$ is 
uniquely determined by $M$. To prove this, we use the fact that $\widehat{R}\otimes_RU$ is an unshrinkable
uniserial $T$-module for each uniserial $R$-module $U$, where $T=\mathrm{End}_R(\widehat{R})$.
When $R$ satisfies a countable condition,
 the collection of uniserial factors of a
  polyserial module $M$ is uniquely determined by $M$
  (Proposition~\ref{P:factor}). 

In this paper all rings are associative and commutative with unity and
all mo\-dules are unital. As in \cite{Fac87} we say that an $R$-module
$E$ is \textbf{divisible} if, for every $r\in
R$ and $x\in E,$ $(0:r)\subseteq (0:x)$ implies that $x\in
rE$, and that $E$ is \textbf{fp-injective}(or
\textbf{absolutely pure}) if $\mathrm{Ext}_R^1(F,E)=0,$ for every finitely
  presented $R$-module $F.$ A ring $R$ is called \textbf{self
    fp-injective} if it is fp-injective as $R$-module. An exact sequence \ $0 \rightarrow F \rightarrow E \rightarrow G \rightarrow 0$ \ is \textbf{pure}
if it remains exact when tensoring it with any $R$-module. In this case
we say that \ $F$ \ is a \textbf{pure} submodule of $E$. Recall that a
module $E$ is fp-injective if and only if it is a pure submodule
of every overmodule. A module is said to be
\textbf{uniserial} if its  submodules are linearly ordered by inclusion
and a ring $R$ is a \textbf{valuation ring} if it is uniserial as
$R$-module. Recall that every finitely presented module over a
valuation ring is a finite direct sum of cyclic modules
\cite[Theorem 1]{War70}. Consequently a module $E$ over a valuation ring
$R$ is fp-injective if and only if it is divisible. 

An $R$-module $F$ is \textbf{pure-injective} if for every pure exact
sequence
\[0\rightarrow N\rightarrow M\rightarrow L\rightarrow 0\]
 of $R$-modules, the following sequence
 \[0\rightarrow\mathrm{Hom}_R(L,F)
\rightarrow\mathrm{Hom}_R(M,F)\rightarrow\mathrm{Hom}_R(N,F)\rightarrow
0\] is exact. An $R$-module $B$
is a \textbf{pure-essential extension} of a submodule $A$ if $A$ is a
pure submodule of $B$ and, if for each
submodule $K$ of $B$, either $K\cap A\ne 0$ or $(A+K)/K$ is not a pure
submodule of $B/K$. 
We say that $B$ is a \textbf{pure-injective hull}
of $A$ if $B$ is pure-injective and a pure-essential extension of $A$.
 By \cite{War69} or \cite[chapter XIII]{FuSa01} each $R$-module $M$ has a
 pure-injective hull and any two pure-injective hulls of $M$ are isomorphic.

In the sequel $R$ is a valuation ring, $P$ its maximal ideal, $Z$ its
subset of zerodivisors and
$\widehat{M}$ the pure-injective hull of $M$, for each $R$-module $M$.
As in \cite[p.69]{FuSa01}, for every proper ideal $A$, we put 
$A^{\sharp}=\{s\in R\mid (A:s)\not=A\}.$ Then $A^{\sharp}/A$ is the set of
zerodivisors of $R/A$ whence $A^{\sharp}$ is a prime ideal. In particular
$\{0\}^{\sharp}=Z$. When $A^{\sharp}=P$, we say that $A$ is an
\textbf{archimedean} ideal. Then $A$ is archimedean if and only if
$R/A$ is self fp-injective.

\section{Properties of $\widehat{R}$} \label{S:R}
The first assertion of the following proposition will play a crucial
role to prove the main results of this paper.

\begin{proposition} \label{P:hull}
 The following assertions hold:
\begin{enumerate}
\item For each $x\in\widehat{R}$
there exist $a\in R,\ p\in P$ and $y\in\widehat{R}$ such that $x=a+pay$.
\item For each archimedean ideal $A$
of $R$, $\widehat{R}/A\widehat{R}$ is an essential extension of $R/A$.
\item $\widehat{R}/P\widehat{R}\cong R/P$.
\end{enumerate}
\end{proposition}
\textbf{Proof.} The third assertion is an immediate consequence of the
first. 

We also deduce the
second assertion from the first. Since $R$ is a pure submodule of
$\widehat{R}$, the natural map 
$R/A\rightarrow \widehat{R}/A\widehat{R}$ is monic. Let
$x\in\widehat{R}\setminus R+A\widehat{R}$. We have
$x=a+pay$ for $a\in R,\ p\in P$ and $y\in\widehat{R}$. Hence $pa\notin
A$. Since $A$ is archimedean, there exists $r\in (A:pa)\setminus
(A:a)$. So $rx\in R+A\widehat{R}\setminus A\widehat{R}$. 

We proceed by steps to prove the first assertion.

\textsl{Step 1.} Suppose that $R$ is self fp-injective. In this case,
$\widehat{R}\cong E_R(R)$ by \cite[Lemma XIII.2.7]{FuSa01}. We may assume that
$x\notin R$. Then there exists $d\in R$ such that $dx\in R$ and $dx\ne
0$. Since $R$ is a pure submodule of $\widehat{R}$ we have $dx=db$ for
some $b\in R$. By \cite[Lemma 2]{Couch03} $(0:x)=(0:b)$, whence $x=bz$ for
some $z\in\widehat{R}$ since $\widehat{R}$ is divisible. In
the same way, there exists $c,\ u\in R$ such that $cz=cu\ne 0$. We get that
$(0:u)=(0:z)=b(0:b)=0$. So $u$ is a unit of $R$. Since $z-u\notin R$,
there exists $s,\ q\in R$ and $y\in\widehat{R}$ such that $0\ne sq=s(z-u)\in R$
and $z-u=qy$. We have $c\in (0:z-u)=(0:q)$. So $q\in P$. Now we put
$a=bu$ and $p=qu^{-1}$ and we get $x=a+pay$.

\textsl{Step 2.} Now we prove that $\widehat{R}/r\widehat{R}\cong
E_{R/rR}(R/rR)$ for each $0\ne r\in P$. If $\cap_{a\ne 0}aR=0$ then it
is an immediate consequence of \cite[Theorem 5.6]{Fac87}. Else $P$ is not
faithful, $R$ is self fp-injective and $\widehat{R}\cong
E_R(R)$. By Step 1 and the implication $1\Rightarrow 2$ the second
assertion holds. So it remains
to show that $\widehat{R}/r\widehat{R}$ is injective over $R/rR$.
Let $J$ be an ideal of $R$ such that $Rr\subset J$ and
$g:J/Rr\rightarrow \widehat{R}/r\widehat{R}$ be a nonzero homomorphism. For
each $x\in\widehat{R}$ we denote by $\bar{x}$ the image of $x$ in
$\widehat{R}/r\widehat{R}$. Let $a\in J\setminus Rr$ such that
$\bar{y}=g(\bar{a})\ne 0$. Then $(Rr:a)\subseteq (r\widehat{R}:y)$. Let
$t\in R$ such that $r=at$. Thus $ty=rz$ for some $z\in\widehat{R}$. It
follows that $t(y-az)=0$. So, since $at=r\ne 0$, we have $(0:a)\subset
Rt\subseteq (0:y-az)$. The injectivity of $\widehat{R}$ implies that there
exists $x\in\widehat{R}$ such that $y=a(x+z)$. We put $x_a=x+z$. If $b\in
J\setminus Ra$ then $a(x_a-x_b)\in r\widehat{R}$. Hence $x_b\in
x_a+(r\widehat{R}:_{\widehat{R}}a)$. Since $\widehat{R}$ is
pure-injective, by \cite[Theorem 4]{War69} there exists $x\in\cap_{a\in
  J}x_a+(r\widehat{R}:_{\widehat{R}}a)$. It follows that
$g(\bar{a})=a\bar{x}$ for each $a\in J$.

\textsl{Step 3.}
Now we prove the first assertion in the general case. If $\cap_{r\ne
   0}rR\ne 0$, then $R$ is self fp-injective. So the result holds by
   Step 1. If
$\cap_{r\ne 0}rR=0$, we put $F=\cap_{r\ne 0}r\widehat{R}$. We will
show that $F=0$. Let $x\in
F\cap R$. Then $x\in R\cap r\widehat{R}=rR$ for each $r\in R,\ r\ne
0$. Therefore $x=0$ and $F\cap R=0$. Let $x\in\widehat{R},\ r,a\in R$ and
$z\in F$ such that $rx=a+z$. There exists $y\in\widehat{R}$ such that
$z=ry$. So $r(x-y)=a$, whence there exists $b\in R$ such that
$rb=a$. It follows that $R$ is a pure submodule of $\widehat{R}/F$. Since
$\widehat{R}$ is a pure-essential extension of $R$ we deduce that
$F=0$. Let $x\in\widehat{R}$. We may assume
that $x\notin R$. There exists $0\ne r\in R$ such that $x\notin
r\widehat{R}$. If $x\in R+r\widehat{R}$ then $x=a+ry$, with $a\in R$ and
$y\in\widehat{R}$. We have $a\notin rR$ else $x\in r\widehat{R}$. So $r=pa$
for some $p\in P$. If $x\notin R+r\widehat{R}$ then, since $R/Rr$ is self fp-injective, from Steps 1 and 2 we
   deduce that $x-a-paz\in r\widehat{R}$ for some $a\in R,\ p\in P$
   and $z\in\widehat{R}$. It is obvious that $a\notin rR$. Now it is
   easy to conclude.  \qed

\bigskip As in the domain case we have:

\begin{proposition} \label{P:flat} $\widehat{R}$ is a faithfully flat
  $R$-module. 
\end{proposition}
\textbf{Proof.} 
Let $x\in\widehat{R}$ and $r\in R$ such that $rx=0$. By
Proposition~\ref{P:hull} there exist $a\in R,\ p\in P$ and
$y\in\widehat{R}$ such that $x=a+pay$. So $rpay\in R$. It follows
that there exists $b\in R$ such that $ra(1+pb)=0$. Hence $ra=0$ and
$r\otimes x=ra\otimes (1+py)=0$. \qed

\section{Pure-injective hulls of uniserial modules} \label{S:unis}
 The following lemma and Proposition~\ref{P:compact} will be
useful to prove the pure-injectivity of some modules in the sequel.
\begin{lemma} \label{L:unis} Let $U$ be a module and $F$ a
  flat module. Then, for each $r,s\in
  R$, $F\otimes_R(sU:_Ur)\cong (F\otimes_RsU:_{F\otimes_RU}r)$.  
\end{lemma}
\textbf{Proof.} We put $E=F\otimes_RU$. Let $\phi$ be the composition of the
  multiplication by $r$ in $U$ with the natural map
$U\rightarrow U/sU$. Then $(sU:_Ur)=\mathrm{ker}(\phi)$. It follows that
  $F\otimes_R(sU:_Ur)$ is isomorphic to
  $\mathrm{ker}(\mathbf{1}_{F}\otimes\phi)$ since $F$ is flat. 
We easily
  check that $\mathbf{1}_{F}\otimes\phi$ is the composition
  of the multiplication by $r$ in $E$ with the natural map
  $E\rightarrow E/sE$. It follows that
$F\otimes_R(sU:_Ur)\cong (sE:_Er)$. \qed

\begin{proposition} \label{P:compact} Every pure-injective $R$-module $F$
  satisfies the following property: if $(x_i)_{i\in I}$ is a family of
  elements of $F$ and $(A_i)_{i\in I}$ a family of ideals of
  $R$ such that the family $\mathcal{F}=(x_i+A_iF)_{i\in I}$
  has the finite intersection property, then $\mathcal{F}$ has a
  non-empty intersection. The converse holds if $F$ is flat.
\end{proposition}
\textbf{Proof.} 
 Let $i\in I$ such that $A_i$ is not finitely
generated. By \cite[Lemma 29]{Couch03} either $A_i=Pr_i$ or
$A_i=\cap_{c\in R\setminus A_i}cR$. If, $\forall i\in I$ such that
 $A_i$ is not finitely generated, we replace  $x_i+A_iF$ by
$x_i+r_iF$ in the first case, and by the family
$(x_i+cF)_{c\in R\setminus A_i}$ in the second case, we deduce from
$\mathcal{F}$ a family $\mathcal{G}$ which has the finite intersection
property. Since $F$ is pure-injective, it follows that there exists
$x\in F$ which belongs to each element of the family
$\mathcal{G}$ by \cite[Theorem 4]{War69}. We may assume that the family
 $(A_i)_{i\in I}$ has no smallest element. So, if $A_i$ is not finitely
 generated, there exists $j\in I$ such that
$A_j\subset A_i$. Let $c\in A_i\setminus PA_j$ such that
 $x_j+cF\in\mathcal{G}$. Then $x-x_j\in cF\subseteq A_iF$ and
 $x_j-x_i\in A_iF$. Hence $x-x_i\in A_iF$ for each $i\in I$. 

Conversely, if $F$ is flat then by Lemma~\ref{L:unis} we have
$(sF:_Fr)=(sR:r)F$ for each $s,r \in R$. We use \cite[Theorem 4]{War69}
to conclude. \qed

\begin{proposition} \label{P:tenspur}
Let $U$ be a uniserial module and $F$ a flat pure-injective
module. Then $F\otimes_RU$ is pure-injective.
\end{proposition}
\textbf{Proof.} Let $E=F\otimes_RU$.
We use \cite[Theorem 4]{War69} to prove
that $E$ is pure-injective. Let
$(x_i)_{i\in I}$ be a family of elements of $F$ such that the family
$\mathcal{F}=(x_i+N_i)_{i\in I}$ has the finite intersection
property, where $N_i=(s_iE:_Er_i)$ and $r_i,s_i\in R$, $\forall i\in I$. 

 First we assume that $U=R/A$ where
$A$ is a proper ideal of $R$. So $E\cong F/AF$. If $s_i\notin
A$ then $N_i=(s_iF:_Fr_i)/AF=(Rs_i:r_i)F/AF$. We set
$A_i=(Rs_i:r_i)$ in this case. If $s_i\in
A$ then $N_i=(AF:_Fr_i)/AF=(A:r_i)F/AF$. We put $A_i=(A:r_i)$ in this
case. For each $i\in I$, let $y_i\in F$ such that $x_i=y_i+AF$. It is obvious that the family $(y_i+A_iF)_{i\in I}$ has the finite
intersection property. By Proposition~\ref{P:compact} this family has a 
non-empty intersection. Then $\mathcal{F}$ has a 
non-empty intersection too. 

Now we assume that $U$ is not finitely generated. It is obvious that
$\mathcal{F}$ has a non-empty intersection if $x_i+N_i=E,\ \forall i\in I$. Now assume there exists $i_0\in I$
such that $x_{i_0}+N_{i_0}\ne E$. Let $I'=\{i\in I\mid N_i\subseteq
N_{i_0}\}$ and $\mathcal{F}'=(x_i+N_i)_{i\in I'}$. Then $\mathcal{F}$ and
  $\mathcal{F}'$ have the same intersection. By
Lemma~\ref{L:unis}
$N_{i_0}=F\otimes_R(s_{i_0}U:_Ur_{i_0})$. It follows that
$(s_{i_0}U:_Ur_{i_0})\subset U$ because $N_{i_0}\ne E$.  Hence  $\exists u\in U$ such that
  $x_{i_0}+N_{i_0}\subseteq F\otimes_RRu$. Then, $\forall
  i\in I'$, $x_i+N_i\subseteq F\otimes_RRu$. We have
  $F\otimes_RRu\cong F/(0:u)F$. From the first part of the proof $F/(0:u)F$ is pure-injective. So we may replace $R$ with
  $R/(0:u)$ and assume that $(0:u)=0$.  Let
  $A_i=((s_iU:_Ur_i):u)$, $\forall i\in I'$. Thus $N_i=A_iF,\ \forall i\in
  I'$. By Proposition~\ref{P:compact}
  $\mathcal{F}'$ has a non-empty intersection. So $\mathcal{F}$ has
  a non-empty intersection too. \qed

\bigskip
Let $U$ be an $R$-module. As in \cite[p.338]{FuSa01} we set 
\[U_{\sharp}=\{s\in R\mid\exists u\in U,\ u\not=0\ \mathrm{and}\
su=0\}\ \mathrm{and}\ U^{\sharp}=\{s\in R\mid sU\subset U\}.\]  
Then $U_{\sharp}$ and $U^{\sharp}$ are prime ideals.

Now it is possible to determine the pure-injective hull of each
uniserial module. We get a generalization of \cite[Corollary XIII.5.5]{FuSa01}

\begin{theorem} \label{T:unis} The following assertions hold:
\begin{enumerate}
\item Let $U$ be a uniserial $R$-module and
  $J=U^{\sharp}\cup U_{\sharp}$. Then
$\widehat{R_J}\otimes_RU$ is the pure-injective hull of
  $U$. Moreover $\widehat{U}$ is an essential
  extension of $U$ if $J=U_{\sharp}$. 
\item For each proper ideal $A$ of $R$,
  $\widehat{R}/A\widehat{R}$ is the pure-injective hull of
  $R/A$. Moreover
  $\widehat{R}/A\widehat{R}\cong E_{R/A}(R/A)$ if $A$ is archimedean.
\end{enumerate}
\end{theorem}
\textbf{Proof.} $(1)$ 
  If $s\in R\setminus J$ then multiplication by $s$ in
$U$ is bijective. So $U$ is an $R_J$-module. After replacing $R$ with
$R_J$, we may assume that $J=P$. We put $\widetilde{U}=\widehat{R_J}\otimes_RU$.

Suppose that $P=U^{\sharp}$. By \cite[Proposition 6]{War69}
$\widetilde{U}=\widehat{U}\oplus V$ where $V$ is a submodule of
$\widetilde{U}$. Let $v\in V$. Then $v=x\otimes u$ where $u\in U$ and
$x\in\widehat{R}$. By Proposition~\ref{P:hull} $x=a+pay$, where $a\in
R$, $p\in P$ and $y\in\widehat{R}$. Since $pU\subset U$, $\exists u'\in
U\setminus (Pu\cup pU)$. Then $u=cu'$ for some $c\in R$ and
$v=cau'+pcay\otimes u'$. We have $y\otimes u'=z+w$ where $w\in V$ and
$z\in\widehat{U}$. So $cau'+pcaz=0$. Since $U$ is pure in
$\widehat{U}$, there exists $z'\in U$ such that $cau'+pcaz'=0$. If
$v\ne 0$ the equality $v=(1+py(\otimes cau'$ implies $cau'\ne 0$. By \cite[Lemma 5]{Couch03} we get that
$u'\in pU$, whence a contradiction. Hence $V=0$.

Now suppose that $P=U_{\sharp}$. If $0\ne z\in \widetilde{U}$ then $z=x\otimes u$
where $u\in U$ and $x\in\widehat{R}$. By
Proposition~\ref{P:hull} there exist $a\in R,\ p\in P$ and
$y\in\widehat{R}$ such that $x=a+pay$. So $z=au+y\otimes pau$. Let
$A=(0:au)$. By \cite[Lemma 26]{Couch03},
$A^{\sharp}=P$. So $(0:pau)=(A:p)\ne A$. Let $r\in (A:p)\setminus
A$. Then $0\ne rz\in U$.

$(2)$ We apply the first assertion by taking $U=R/A$. In this case,
$U^{\sharp}=P$. The pure-injective hull of $R/A$ is the same over $R$
and over $R/A$. Since $R/A$ is self fp-injective when $A$ is archimedean
then we use \cite[Lemma XIII.2.7]{FuSa01} to prove the last assertion.
\qed

\bigskip In the previous theorem, if $U$ is not cyclic and if
$U^{\sharp}\subseteq U_{\sharp}$ then $\widehat{U}$ is not necessarily
isomorphic to $E_{R/(0:U)}(U)$. For instance:

\begin{example} \textnormal{Assume that $P=Z$ and $P$ is faithful. We
    choose $U=P$. Then $U^{\sharp}=U_{\sharp}=P$,
    $\widehat{U}=P\widehat{R}$ and $E_R(U)=\widehat{R}$.}
\end{example}

If $U$ is a non-standard uniserial module over a valuation domain $R$ then
$\widehat{U}$ is indecomposable by \cite[Proposition 5.1]{Fac87} and
 there exists a standard uniserial module $V$ such that $\widehat{U}\cong\widehat{V}$ 
by \cite[Theorem XIII.5.9]{FuSa01}.
So, $\widehat{R}\otimes_RU\cong\widehat{R}\otimes_RV$ doesn't imply $U\cong V$. 
However, it is possible to get the following proposition:
\begin{proposition} \label{P:isounis}
Let $U$ and $V$ be uniserial modules and $J=U^{\sharp}\cup U_{\sharp}$. Assume that
$\widehat{R}\otimes_RU\cong\widehat{R}\otimes_RV$. Then $U$ and $V$ are isomorphic if one
of the following conditions is satisfied:
\begin{enumerate}
\item $U^{\sharp}=J$ and $J\ne J^2$,
\item $U$ is countably generated.
\end{enumerate}
\end{proposition}
\textbf{Proof.} Let $\phi:\widehat{R}\otimes_RU\rightarrow\widehat{R}\otimes_RV$ be the
 isomorphism. Let $0\ne u\in U$. Then $\phi(u)=x\otimes v$ for some $x\in\widehat{R}$ and
 $v\in V$. By proposition~\ref{P:hull} we may assume that $x=1+py$ for some $p\in P$ and
$y\in\widehat{R}$. First we shall prove that $(0:u)=(0:v)$. It is obvious that 
$(0:v)\subseteq (0:u)$. Let $r\in (0:u)$. Then $x\otimes rv=0$. From the flatness of 
$\widehat{R}$ we deduce that there exist $s\in R$ and $z\in\widehat{R}$ such that $x=sz$ and 
$srv=0$. If $s\in P$ then we get that $1=qe$ for some $q\in P$ and $e\in\widehat{R}$. Since
 $R$ is pure in $\widehat{R}$, it follows that $1\in P$. This is absurb. Hence $s$ is a unit 
and  $r\in (0:v)$.

Let $v,v'$ be nonzero elements of $ V$ and $x,x'\in 1+P\widehat{R}$  such
 that $x\otimes v=x'\otimes v'$. There exists $t\in R$ such that  $v=tv'$. Now we shall prove that
$t$ is a unit of $R$. We get that $(x'-tx)\otimes v'=0$. If $t\in P$, as above we deduce that
$v'=0$, whence a contradiction.

Let $u\in U$ and $v\in V$ as in the first part of the proof. By \cite[Lemma 26]{Couch03} we have
$U_{\sharp}=(0:u)^{\sharp}=(0:v)^{\sharp}=V_{\sharp}$. Let $p\in P$. We shall prove that
 $u\in pU$ if and only if $v\in pV$. If $v=pw$ for some $w\in V$ then
 $\phi(u)=px\otimes w=p\phi(z)$ for some
 $z\in\widehat{R}\otimes_RU$. Since $U$ is a pure submodule, then $u=pu'$ for some $u'\in U$.
Conversely, if $u=pu'$ for some $u'\in U$ and $\phi(u')=x'\otimes v'$ where $v'\in V$ and
$x'\in 1+P\widehat{R}$, we get that $x'\otimes pv'=x\otimes v$. From above, we deduce that
$v\in pV$. So, $U^{\sharp}=V^{\sharp}$.

Now we can prove that $U$ and $V$ are isomorphic when the first condition is satisfied. In
 this case $U$ and $V$ are modules over $R_J$. Since $J\ne J^2$ , $JR_J$ is a 
principal ideal of $R_J$. Since $JU\subset U$ and $JV\subset V$, $U$ and $V$ are cyclic over
$R_J$. Let $u\in U$ and $v\in V$ as in the first part of the proof, and suppose that $U=R_Ju$.
If $v=rw$ for some $r\in R_J$ and $w\in V$ then we get, as above, that $u=ru'$ for some $u'\in U$. 
So $r$ is a unit and $U$ and $V$ are isomorphic.

Let $\{u_i\}_{i\in I}$ be a spanning set of $U$. For each $i\in I$, let $v_i\in V$ 
and $x_i\in 1+P\widehat{R}$ such that $\phi(u_i)=x_i\otimes v_i$. Suppose that
 $(0:U)\subset (0:u),\ \forall u\in U$. From the first part of
proof we deduce that $(0:V)\subset (0:v),\ \forall v\in V$. We have
 $\cap_{i\in I}(0:u_i)=(0:U)$. Thus $\cap_{i\in I}(0:v_i)=(0:V)$. So, for each $v\in V$
there exists $i\in I$ such that $(0:v_i)\subset (0:v)$. Hence $v\in Rv_i$. Now, suppose 
$\exists u\in U$ such that $(0:u)=(0:U)$. By \cite[Lemma X.1.4]{FuSa01} $J=U^{\sharp}$. We may 
assume that $J=J^2$ and $I$ is infinite. Then $JU=U$ and $JV=V$. Let $v\in V$. There exists
$p\in J$ such that $v\in pV$. But there exists $i\in I$ such that $u_i\notin pU$. So,
$v_i\notin pV$. Hence $v\in R_Jv_i$. Now suppose that $I=\mathbb{N}$. Let 
$(a_n)_{n\in\mathbb{N}}$ be a sequence of elements of $P$ such that 
$u_n=a_nu_{n+1},\ \forall n\in\mathbb{N}$. We put $\varphi(u_0)=v_0$. Suppose that 
$\varphi(u_n)=s_nv_n$ where $s_n$ is a unit. By the second part of the proof there exists a
 unit $t_n$ such that $a_nv_{n+1}=t_n\varphi(u_n)$. Hence we set 
$\varphi(u_{n+1})=t_n^{-1}v_{n+1}$. So, by induction on $n$,  we get
an isomorphism $\varphi:U\rightarrow V$. \qed

\bigskip Let $T=\mathrm{End}_R(\widehat{R})$. Then $T$ is a local ring by \cite[Proposition 5.1]{Fac87} and 
\cite[Theorem XIII.3.10]{FuSa01}. For each $R$-module $M$,
$\widehat{R}\otimes_RM$ is a left $T$-module. As in \cite{Fac98} we say that a left uniserial $T$-module $F$
is \textbf{shrinkable} if there exists two $T$-submodules $G$ and $H$ of $F$ such that $0\subset H\subset G\subset F$
and $F\cong G/H$. Otherwise $F$ is said to be \textbf{unshrinkable}. 
\begin{proposition} \label{P:end}
Let $U$ be a uniserial $R$-module. Then:
\begin{enumerate}
\item $\widehat{R}\otimes_RU$ is a left unshrinkable uniserial $T$-module.
\item $\mathrm{End}_T(\widehat{R}\otimes_RU)$ is a local ring.
\end{enumerate}
\end{proposition}
\textbf{Proof.} $(1)$ Let $x\in 1+P\widehat{R}$. First we prove that $Rx$ is a pure submodule of $\widehat{R}$.
Let $a,b\in R$ and $y\in\widehat{R}$ such that $by=ax$. By Proposition~\ref{P:hull} $y=c+pcz$ for some $c\in R$,
$p\in P$ and $z\in\widehat{R}$. Suppose that $a\notin Rbc$. Then $bc=ra$ for some $r\in P$. If $x=1+qx'$ for some
$q\in P$ and $x'\in\widehat{R}$, we get that $a(1-r)=a(rpz-qx')=aty'$ for some $t\in P$ and $y'\in\widehat{R}$.
Since $R$ is a pure submodule of $\widehat{R}$ there exists $s\in R$ such that $a(1-r-ts)=0$. We deduce that $a=0$,
whence a contradiction. So $a\in Rbc$. By using similar arguments we easily show that $Rx$ is faithful.

Let $z,z'\in\widehat{R}\otimes_RU$. We have $z=x\otimes u$ and $z'=x'\otimes u'$ where $x,x'\in 1+P\widehat{R}$ and
$u,u'\in U$. Assume that $u'=ru$ for some $r\in R$. The homomorphism $\phi:Rx\rightarrow Rrx'$ such that $\phi(x)=rx'$ is well defined
and can be extended to $\widehat{R}$. We get that $\phi z=z'$. Hence $\widehat{R}\otimes_RU$ is uniserial over $T$.

Suppose that $\widehat{R}\otimes_RU$ is shrinkable over $T$. By \cite[Lemma 1.17]{Fac98} there exists $z\in\widehat{R}\otimes_RU$
such that $Tz$ is shrinkable. We have $z=x\otimes u$ where $x\in 1+P\widehat{R}$ and $u\in U$. So $Tz=\widehat{R}\otimes_RRu$.
There exist $z'\in Tz$ and a non-injective $T$-epimorphism $\alpha:Tz'\rightarrow Tz$. Let $K=\mathrm{Ker}\ \alpha$. We may assume
that $\alpha(z')=z$. We have $z'=x'\otimes ru$ where $x'\in 1+P\widehat{R}$ and $r\in R$. Let $y$ be a nonzero element of $K$. Thus 
$y=tz'=ay'\otimes ru$ for some $t\in T$, $y'\in 1+P\widehat{R}$ and $a\in R$. But there exist $s,s'\in T$ such that $x'=sy'$ and $y'=s'x'$.
So $0\ne ax'\otimes ru\in K$. Since $y\ne 0$ we have $aru\ne 0$. On the other hand $x\otimes aru=\alpha(ax'\otimes ru)=0$. It follows that
$aru=0$ whence a contradiction. So $\widehat{R}\otimes_RU$ is unshrinkable.

$(2)$ is an immediate of $(1)$ and \cite[Proposition 9.24]{Fac98}.
\qed

\begin{proposition} \label{P:gene}
Let $\mathfrak{c}$ be a cardinal. Consider a $\mathfrak{c}$-generated $R$-module $M$ and
$U$ a pure uniserial $R$-submodule of $M$. Then $U$ is $\mathfrak{c}$-generated.
\end{proposition}
\textbf{Proof.} We easily check that $\widehat{R}\otimes_RU$ is a pure submodule of $\widehat{R}\otimes_RM$.
By Proposition~\ref{P:tenspur} $\widehat{R}\otimes_RU$ is pure-injective. Hence  $\widehat{R}\otimes_RU$ is a summand of $\widehat{R}\otimes_RM$.
On the other hand $\widehat{R}\otimes_RM$ is a $\mathfrak{c}$-generated $T$-module. Then $\widehat{R}\otimes_RU$
is also $\mathfrak{c}$-generated over $T$. We may assume that $\widehat{R}\otimes_RU$ is generated by $(1\otimes u_i)_{i\in I}$, 
where $I$ is a set whose cardinal is $\mathfrak{c}$ and $u_i\in U,\ \forall i\in I$. Let $V$ be the submodule of $U$ generated by
$(u_i)_{i\in I}$. Then the inclusion map $V\rightarrow U$ induces an isomorphism $\widehat{R}\otimes_RV\rightarrow\widehat{R}\otimes_RU$.
Since $\widehat{R}$ is faithfully flat it follows that $V=U$. \qed

\bigskip From Theorem~\ref{T:unis} we deduce the following corollary
on the structure of indecomposable injective modules.
\begin{corollary} Let $E$ be an indecomposable injective module,
  $J=E_{\sharp}$ and $\mathcal{A}(E)=\{(0:_{R_J}x)\mid 0\ne
  x\in E\}$. Then:
\begin{enumerate}
\item $\forall A,B\in\mathcal{A}(E),\ A\subseteq B$ there exists a
  monomorphism 
\[\varphi_{A,B}:\widehat{R_J}/B\widehat{R_J}\rightarrow
  \widehat{R_J}/A\widehat{R_J}\] 
such that
  $\varphi_{A,C}=\varphi_{A,B}\circ\varphi_{B,C}$, $\forall
  A,B,C\in\mathcal{A}(E),\ A\subseteq B\subseteq C$.
\item
  $E\cong\varinjlim\{(\widehat{R_J}/A\widehat{R_J},\varphi_{A,B})\mid
  A,B\in\mathcal{A}(E),\ A\subseteq B\}$.
\item  $E\cong\widehat{R_J}/(0:_{R_J}e)\widehat{R_J}$ if
  $(0:_{R_J}e)=(0:_{R_J}E)$ for some $e\in E$.
\item Suppose that $E$ contains a uniserial
    $R_J$-module $U$ such that $\mathcal{A}(E)=\mathcal{A}(U)$\footnote{We
    know that this condition holds if $R$ satisfies an additional
    hypothesis: see \cite[Corollary 22]{Couch03}, \cite[Theorem
    5.5]{ShLe74} or Remark~\ref{R:countable}}. Then $E\cong\widehat{R_J}\otimes_RU$. Moreover,
    $\forall A,B\in\mathcal{A}(E),\ A\subseteq B$, there exists
    $r\in R$ such that one can
    choose $\varphi_{A,B}=\mathbf{1}_{\widehat{R_J}}\otimes\bar{r}$ where
    $\bar{r}:R_J/B\rightarrow R_J/A$ is defined by
    $\bar{r}(a+B)=ar+A,\ \forall a\in R$.
\end{enumerate}
\end{corollary}
\textbf{Proof.} $(1)$ If $A\in\mathcal{A}(E)$ then $A^{\sharp}=J$ by
\cite[Lemma 26]{Couch03}. So $A$ is an archimedean ideal
of $R_J$. By Theorem~\ref{T:unis} there exists an isomorphism
\[\phi_A:\widehat{R_J}/A\widehat{R_J}\rightarrow (0:_EA).\]
 Let
$u_{A,B}:(0:_EB)\rightarrow (0:_EA)$ be the inclusion map , $\forall A,B\in\mathcal{A}(E),\ A\subseteq B$. We set
$\varphi_{A,B}=\phi_A^{-1}\circ u_{A,B}\circ\phi_B$. It is easy to
check the first assertion. 

$(2)$ and $(3)$  These assertions are now obvious. 

$(4)$ First we prove that $U$ is fp-injective. Let $x\in E$ and $s\in
R$ such that $0\ne sx\in U$. We put $u=sx$. From
$\mathcal{A}(E)=\mathcal{A}(U)$, it follows that $\exists v\in U$ such
that $(0:_{R_J}v)=(0:_{R_J}x)$ and consequently $u=tv$ for some $t\in
R$. We set $A=(0:_{R_J}x)$. We get that
$(0:_{R_J}u)=(A:_{R_J}t)=(A:_{R_J}s)$. By \cite[Lemma 26]{Couch03}
$A^{\sharp}=E_{\sharp}=J$. It follows that $R_Js=R_Jt$. So $U$ is a
pure submodule of $E$. We conclude by Theorem~\ref{T:unis} and
\cite[Lemma XIII.2.7]{FuSa01} that $E\cong \widehat{R_J}\otimes_RU$.

Let $u,v\in
    U$ such that $(0:_{R_J}u)=A$ and $(0:_{R_J}v)=B$. There exists
    $r\in R$ such that  $v=ru$ and $B=(A:r)$ (if $A=B$ we take $v=u$
    and $r=1$). So $\bar{r}$ is a monomorphism.
\qed

\section{Pure-injective hulls of polyserial modules} \label{S:poly}
We say that a module $M$ is \textbf{polyserial} if it has a
pure-composition series
\[0=M_0\subset M_1\subset\dots\subset M_n=M,\]
(i.e. $M_k$ is a pure submodule of $M$, for each $k$, $0\leq k\leq n$)
where $M_k/M_{k-1}$ is uniserial for each $k$, $1\leq k\leq n$. By
\cite[Lemma I.7.8]{FuSa01}, if $M$ is finitely generated, $M$ has a
pure-composition series, where $M_k/M_{k-1}\cong R/A_k$ and $A_k$ is a proper ideal, for each
$k,\ 1\leq k\leq n$. We denote by $\mathrm{gen}\ M$ the
minimal number of generators of $M$. By \cite[Lemma V.5.3]{FuSa01}
$n=\mathrm{gen}\ M$. The following sequence $(A_1,\cdots,A_n)$ is
called the \textbf{annihilator sequence} of $M$ and is uniquely
determined by $M$, up to the order (see \cite[Theorem V.5.5]{FuSa01}). 

 Now we can extend the result obtained by Warfield\cite{War69} in
the domain case for finitely generated modules.
\begin{theorem} \label{T:main} Let $M$ be a finitely generated
  $R$-module. Then $\widehat{R}\otimes_RM\cong\widehat{M}$. Moreover,
  $\widehat{M}\cong
  \widehat{R}/A_1\widehat{R}\oplus\cdots\oplus\widehat{R}/A_n\widehat{R}$ where
  $(A_1,\cdots,A_n)$ is the annihilator sequence of $M$.
\end{theorem}
\textbf{Proof.} It is easy to verify that $M$ is a pure submodule of
$\widehat{R}\otimes_RM$.  We have that
$\widehat{R}\otimes_RM_1$ is a pure submodule of
$\widehat{R}\otimes_RM$ too. By
Proposition~\ref{P:tenspur} $\widehat{R}\otimes_RM_1$ is pure-injective. It
follows that
$\widehat{R}\otimes_RM\cong(\widehat{R}\otimes_RM_1)\oplus(\widehat{R}\otimes_RM/M_1)$.
By induction on $n$ we get that $\widehat{R}\otimes_RM\cong
  \widehat{R}/A_1\widehat{R}\oplus\cdots\oplus\widehat{R}/A_n\widehat{R}$. So
  $\widehat{R}\otimes_RM$ is pure-injective. By \cite[Proposition 6]{War69}
  $\widehat{M}$ is a direct summand of $\widehat{R}\otimes_RM$. So
  $\widehat{R}\otimes_RM\cong\widehat{M}\oplus V$, where $V$ is a submodule
  of $\widehat{R}\otimes_RM$. From Proposition~\ref{P:hull} we deduce that,
  for each
  $x\in\widehat{R}\otimes_RM$, there exist $m\in M,\ p\in P$ and
  $y\in\widehat{R}\otimes_RM$ such that $x=m+py$. Assume that $x\in
  V$. There exists $z\in\widehat{M}$ and $v\in V$ such that
  $x=m+pz+pv$. It follows that $x=pv$, whence $V=PV$. On the other
  hand, $\widehat{R}/A\widehat{R}$ is indecomposable by \cite[Proposition
  5.1]{Fac87} and $\mathrm{End}_R(\widehat{R}/A\widehat{R})$ is local by
  \cite[Theorem 9]{ZHZ78} or \cite[Theorem XIII.3.10]{FuSa01}, for every
  proper ideal $A$. By Krull-Schmidt Theorem 
  $V\cong\widehat{R}/A_{k_1}\widehat{R}\oplus\cdots\oplus\widehat{R}/A_{k_p}\widehat{R}$
 where $\{k_1,\cdots,k_p\}$ is a subset of $\{1,\cdots,n\}$. If $V\ne 0$, by Proposition~\ref{P:hull}
 we get $V\ne PV$. This contradiction completes the proof. \qed

\bigskip 
The \textbf{Malcev rank} of a module $N$ is defined as the cardinal
number
\[\mathrm{Mr}\ N=\mathrm{sup}\{\mathrm{gen}\ M\mid M \subseteq N,\ 
\mathrm{gen}\ M<\infty\}.\] 
The following proposition is identical to the first part of \cite[Proposition
XII.1.6]{FuSa01}. Here we give a different proof.
\begin{proposition} \label{P:length} The length of any
  pure-composition series of a polyserial module $M$ equals
  $\mathrm{Mr}\ M$.
\end{proposition}
\textbf{Proof.} Let $0=M_0\subset M_1\subset\dots\subset M_n=M$ be a
pure-composition series of $M$ with uniserial factors. As in \cite[Corollary
XII.1.5]{FuSa01} we prove that $\mathrm{Mr}\ M\leq n$. Equality holds for $n=1$. From
the pure-composition series of $M$, we deduce a pure-composition
series of $M/M_1$ of length $n-1$. By induction hypothesis $M/M_1$ contains a
finitely generated submodule $Y$ with $\mathrm{gen}\ Y=n-1$.

Assume that $Y$ is generated by
$\{y_2,\ldots,y_n\}$. Let
$x_2,\ldots,x_n\in M$  such that $y_k = x_k+M_1$ and $F$ be the
submodule of $M$ generated by $x_2,\ldots,x_n$. 
If $F\cap M_1=M_1$ then  $M_1\subseteq F$ and $M_1$ is a pure
submodule of $F$. In this case $M_1$ is finitely generated by Proposition~\ref{P:gene}.
 It follows that the following sequence
is exact:
\[0\rightarrow
\frac{M_1}{PM_1}\rightarrow\frac{F}{PF}\rightarrow\frac{Y}{PY}\rightarrow 0.\] 
So we have $\mathrm{gen}\ Y=\mathrm{gen}\ F-\mathrm{gen}\ M_1\leq
n-2$. We get a contradiction since $\mathrm{gen}\ Y=n-1$. Hence $F\cap
M_1\not= M_1$. Let $x_1\in M_1\setminus F\cap M_1$. Let $X$ be the submodule of
$M$ generated by $x_1,\ldots,x_n$. Clearly $Rx_1=M_1\cap X$. We will show that $Px_1=Rx_1\cap
PX$. Let $x\in Rx_1\cap PX$. Then $x=p\sum_{k=1}^{k=n}a_kx_k=rx_1$
where $p\in P$ and $r,\ a_1,\dots,\ a_n$ are elements of $R$. It
follows that $p\sum_{k=2}^{k=n}a_kx_k=(r-pa_1)x_1$. So $(r-pa_1)x_1\in
M_1\cap F\subset Rx_1$. We deduce that $r-pa_1\in P$ whence $r\in
P$. Hence $x\in Px_1$. Consequently the following sequence is exact:
\[0\rightarrow
\frac{Rx_1}{Px_1}\rightarrow\frac{X}{PX}\rightarrow\frac{Y}{PY}\rightarrow 0.\]
Then $\mathrm{gen}\ X=n$. \qed

\bigskip Now we study the pure-injective hulls of polyserial modules. 

\begin{theorem} \label{T:poly} Let $M$ be a polyserial module with the
    following pure-composition series:
\[0=M_0\subset M_1\subset\dots\subset M_n=M\]
For each integer $k$, $1\leq k\leq n$ we put $U_k=M_k/M_{k-1}$. Then:
\begin{enumerate}
\item There exists a subset $I$ of $\{k\in\mathbb{N}\mid 1\leq k\leq n\}$
such that $\widehat{M}\cong\oplus_{k\in I}\widehat{U_k}$.
\item
  $\widehat{R}\otimes_RM$ is pure-injective and isomorphic to
  $\oplus_{k=1}^{k=n}\widehat{R}\otimes_RU_k$. 
\item The collection $(\widehat{R}\otimes_RU_k)_{1\leq k\leq n}$ is uniquely determined by $M$.

\end{enumerate}
\end{theorem}
\textbf{Proof.} $(1)$ Let $N$ be a pure submodule of $M$. The inclusion map
$N\rightarrow \widehat{N}$ can be extended to $w:M\rightarrow
\widehat{N}.$ Let $f:M\rightarrow \widehat{N}\oplus \widehat{M/N}$ defined by
$f(x)=(w(x),x+N),$ for each $x\in M.$ It is easy to verify that
$f$ is a pure monomorphism. It follows that $\widehat{M}$ is a summand
of $\widehat{N}\oplus\widehat{M/N}$. So, by induction on $n$, we
easily get that $\widehat{M}$ is a summand of
$\oplus_{k=1}^{k=n}\widehat{U_k}$. Since, $\forall k\in\mathbb{N}$, $1\leq k\leq
n$, $\widehat{U_k}$ is
indecomposable by \cite[Proposition 5.1]{Fac87} and
$\mathrm{End}_R(\widehat{U_k})$ is local by \cite[Theorem 9]{ZHZ78} or
\cite[Theorem XIII.3.10]{FuSa01}, we apply Krull-Schmidt Theorem to conclude. 

$(2)$ We do as in the proof of Theorem~\ref{T:main}. 

$(3)$ Since $\widehat{R}\otimes_RM$ and $\widehat{R}\otimes_RU_k$ are $T$-modules,
 we conclude by Proposition~\ref{P:end} and Krull-Schmidt theorem.\qed

\begin{corollary} \label{C:counpoly}
If $M$ is polyserial and countably generated
then any two pure-composition series of $M$ are isomorphic.
\end{corollary}
\textbf{Proof.} By Theorem~\ref{T:poly}
the collection
    $(\widehat{R}\otimes_RU_k)_{1\leq k\leq n}$ is uniquely determined
    by $M$. It remains to show that, if $U$ and $V$ are uniserial
    modules such that $\widehat{R}\otimes_RU\cong \widehat{R}\otimes_RV$
    then $U\cong V$.  It is an immediate consequence of
    Proposition~\ref{P:gene} and Proposition~\ref{P:isounis}.
     \qed

\bigskip
Recall that an $R$-module $M$ is \textbf{finitely} (respectively
\textbf{countably cogenerated}) if $M$ is a submodule of a product
 of finitely (respectively countably) many injective hulls of simple modules.

 The following proposition completes \cite[Corollary 35]{Couch03}.
\begin{proposition} \label{P:countable}
The following conditions are equivalent:
\begin{enumerate}
\item Every finitely generated $R$-module is countably cogenerated and
  every ideal of $R$ is countably generated.
\item For each prime ideal $J$ which is the union of the
  set of primes properly contained in $J$ there is a countable subset
  whose union is $J,$ and for each prime ideal $J$ which
  is the intersection of the
  set of primes containing properly $J$ there is a countable subset
  whose intersection is $J.$
\item Each uniserial module is countably generated.
\end{enumerate}
\end{proposition}
$(1)\Leftrightarrow (2)$ holds by \cite[Corollary 35]{Couch03}

$(3)\Rightarrow (2)$ Let $J$ be a prime ideal. Then $J$ and $R_J$ are
uniserial $R$-modules. So they are countably generated. If $R_J$ is
generated by $\{t_n^{-1}\mid n\in\mathbb{N}\}$, where $t_n\notin J$
$\forall n\in\mathbb{N}$, then $J=\cap_{n\in\mathbb{N}}Rt_n$. Now it is
easy to get the second condition.

$(1)\Rightarrow (3)$ Let $U$ be a uniserial module and $J=U^{\sharp}\cup
U_{\sharp}$. Then $U$ is an $R_J$-module. But $R/J$ countably
cogenerated is equivalent to $R_J$ countably generated. Hence $U$ is
countably generated over $R$ if and only if $U$ is countably generated
over $R_J$. So we may assume that $J=P$.

First assume that $U^{\sharp}=P$. If $PU\subset U$ then $U=Ru$ where
$u\in U\setminus PU$. Now suppose that $PU=U$. Let $r,s\in P$ such
that $rU\ne 0$. If $rU=rsU$ then by \cite[Lemma 5]{Couch03} we have
$U=sU$, hence a contradiction. Let $\{p_n\mid n\in\mathbb{N}\}$ be a
spanning set of $P$ such that $p_{n+1}\notin Rp_n$. Then
$U=\cup_{n\in\mathbb{N}}p_nU$. We may assume that $p_nU\ne 0,\ \forall
n\in\mathbb{N}$. So $p_nU\subset p_{n+1}U$ for each 
$n\in\mathbb{N}$. Let $u_n\in p_{n+1}U\setminus p_nU$ for each
$n\in\mathbb{N}$. Then $U$ is generated by $\{u_n\mid
n\in\mathbb{N}\}$.

Now suppose that $U_{\sharp}=P$. Assume that $(0:u)=(0:U)$ for some
$u\in U$. Let $v\in U$ such that $u=av$ for some $a\in R$. By
\cite[Lemma 2]{Couch03} $(0:u)=((0:v):a)$. We get that
$(0:v)=((0:v):a)=(0:U)$. Since $(0:v)^{\sharp}=P$ by \cite[Lemma
26]{Couch03} $a$ is a unit, and consequently $U$ is cyclic. Now we assume
that $(0:U)\subset (0:u)$ for each $u\in U$. We have $(0:U)=\cap_{u\in
  U}(0:u)$. By \cite[Lemma 30]{Couch03} there exists a countable family
$(u_n)_{n\in\mathbb{N}}$ of elements of $U$ such that
$(0:U)=\cap_{n\in\mathbb{N}}(0:u_n)$ and $u_{n+1}\notin Ru_n$, $\forall
n\in\mathbb{N}$. If $u\in U$, since  $(0:u)\ne (0:U)$, then there
exists $n\in\mathbb{N}$ such that $(0:u_n)\subset (0:u)$. Hence $u\in
Ru_n$ and $U$ is generated by $\{u_n\mid n\in\mathbb{N}\}$.
\qed

\begin{remark} \label{R:countable}
\textnormal{In the same way, one can prove that the two first
 conditions of \cite[Proposition 32]{Couch03} (respectively
 \cite[Corollary 34]{Couch03}) are equivalent to the following: each
 indecomposable injective module $E$ such that $E_{\sharp}=P$ contains a
 uniserial pure submodule which is countably generated (respectively  each
 indecomposable injective module contains a uniserial pure submodule which is
 countably generated).}
\end{remark}

\begin{proposition} \label{P:factor}
Suppose that $R$ satisfies the equivalent conditions of
  Proposition~\ref{P:countable}. Then any two pure-composition series
  of a polyserial $R$-module are isomorphic.
\end{proposition}
\textbf{Proof.}  It is an immediate consequence of Proposition~\ref{P:countable}
and Corollary~\ref{C:counpoly}.
     \qed

\section{Two criteria for maximality of $R$} \label{S:max}
\bigskip By Theorem~\ref{T:main}, if $M$ is finitely generated,
then $\widehat{M}$ is a direct sum of $\mathrm{gen}\ M$ indecomposable
pure-injective modules and $\mathrm{gen}\ M=\mathrm{Mr}\ M$ by
\cite[Corollary XII.1.7]{FuSa01}. But Theorem~\ref{T:malcev} proves
that, if $M$ is polyserial, then $\widehat{M}$ is not necessarily a
direct sum of $\mathrm{Mr}\ M$ indecomposable pure-injective modules.

 As in \cite{SaZa85}, if $x\in\widehat{R}\setminus R$,
 we say that $\mathrm{B}(x)=\{r\in R\mid x\notin R+r\widehat{R}\}$ is the
 \textbf{breath ideal} of $x$. Then
 Proposition~\ref{P:breath} is a generalization of 
\cite[Proposition 1.4]{SaZa85}. The following lemma is useful to prove
 this proposition. 

\begin{lemma} \label{L:inters} 
Let $J$ be a proper ideal such that $J=\cap_{c\notin
  J}cR$. Then $J\widehat{R}=\cap_{c\notin J}c\widehat{R}$.
\end{lemma}
\textbf{Proof.}
 By Theorem~\ref{T:unis} $\widehat{R}/J\widehat{R}$ is the
 pure-injective hull of $R/J$. In the proof of Step 3 of
 Proposition~\ref{P:hull} it is already shown that $\cap_{a\ne
 0}a\widehat{R}=0$ if $\cap_{a\ne 0}aR=0$. So we apply this result to
 $R/J$ to get the lemma.
\qed

\bigskip Recall that the \textbf{ideal topology} of $R$ is the linear
topology which has as a basis of neighborhoods of $0$ the nonzero
principal ideals.
\begin{proposition} \label{P:breath} Let $A$ be a proper ideal. Then
 $R/A$ is Hausdorff and non-complete in its ideal topology if and only
 if $A=\mathrm{B}(x)$ for some $x$ in $\widehat{R}\setminus R$.
\end{proposition}
\textbf{Proof.} To show that $R/\mathrm{B}(x)$ is Hausdorff, we do as
in \cite[Proposition 1.4]{SaZa85}, we prove that $a\notin\mathrm{B}(x)$
implies that $pa\notin\mathrm{B}(x)$ for some $p\in P$. We have
$x=r+ay$ where $r\in R$ and $y\in\widehat{R}$. By Proposition~\ref{P:hull},
$\widehat{R}=R+P\widehat{R}$. So $y=s+pz$, for some $s\in R$, $p\in P$
and $z\in\widehat{R}$. Therefore we get $x=r+as+paz\in
R+pa\widehat{R}$. For each $a\notin\mathrm{B}(x)$, $x\in
r_a+a\widehat{R}$ for some $r_a\in R$. If the family
$(r_a+aR)_{a\notin\mathrm{B}(x)}$ has a non-empty intersection then,
by using Lemma~\ref{L:inters}, we get that
  $x\in R+\mathrm{B}(x)\widehat{R}$, whence a contradiction. So
  $R/\mathrm{B}(x)$ is non-complete.

Conversely, assume that $R/A$ is Hausdorff and non-complete. Then
there exists a family $(r_a+aR)_{a\notin A}$ which has the finite
intersection and an empty total intersection. Since $\widehat{R}$ is
pure-injective, the total intersection of the family
$(r_a+a\widehat{R})_{a\notin A}$ contains an element $x$ which doesn't
belong to $R$. Clearly $\mathrm{B}(x)\subseteq A$. If
$x=r+b\widehat{R}$ for some $r\in R$ and $b\in A$ then $r\in r_a+aR$
for each $a\notin A$, since $R$ is a pure submodule of
$\widehat{R}$. We get a contradiction. So $A=\mathrm{B}(x)$. 
\qed

\bigskip The following lemma is a generalization of \cite[Lemma
1.3]{SaZa85}. It will be useful to prove Theorem~\ref{T:maximal}.
\begin{lemma} \label{L:breath} Let $x\in\widehat{R}$ such that
  $x=r+ay$ for some $r,a\in R$ and $y\in\widehat{R}$. Then
  $\mathrm{B}(y)=(\mathrm{B}(x):a)$.
\end{lemma}
\textbf{Proof.} Let $t\notin\mathrm{B}(y)$. Then $y=s+tz$ for some $s\in
R$ and $z\in\widehat{R}$. It follows that $x=r+as+aty$. So $t\notin
(\mathrm{B}(x):a)$. 

Conversely, if $t\notin (\mathrm{B}(x):a)$ then we get the following
equalities $x=r+ay=s+taz$ for some $s\in R$ and
$z\in\widehat{R}$. Since $R$ is a pure submodule of $\widehat{R}$ it
follows that $a(y-tz-b)=0$ for some $b\in R$. From the flatness of
$\widehat{R}$ we deduce that $(y-tz-b)\in (0:a)\widehat{R}$. But
$ta\notin\mathrm{B}(x)$ implies that $ta\ne 0$, whence
$(0:a)\subset Rt$. Hence $t\notin\mathrm{B}(y)$. \qed

\begin{theorem} \label{T:maximal} Let $N$ be the nilradical of
  $R$. Then $R$ is maximal if and only if $R/N$ and $R_N$ are maximal.
\end{theorem}
\textbf{Proof.} Suppose that $R$ is maximal. It is obvious that $R/N$ is
maximal. By \cite[Lemma 2]{Gil71} $R_N$ is maximal too.

Conversely assume that $R/N$ and $R_N$ are maximal. Let $K$ be the
kernel of the natural map $R\rightarrow R_N$. Let $r\in K$. Thus there exists $s\in R\setminus N$ such that $sr=0$. It follows that $K\subseteq N\subset (0:r)$. Then $K^2=0$. So $K$ is
a uniserial $R/K$-module which is linearly compact if $R/K$ is maximal. Consequently $R$ is maximal if and only if $R/K$ is
maximal. In the sequel we may assume that $K=0$. So $N=Z$ and it is an
$R_N$-module. It is enough to show that $N$ is a linearly compact
module. Let $(A_i)_{i\in I}$ be a family of ideals contained in $N$ and
$(x_i)_{i\in I}$ a family of elements of $N$ such that the family
$\mathcal{F}=(x_i+A_i)_{i\in I}$ has the finite intersection
property. We put $A=\cap_{i\in I}A_i$. We may assume that $A\subset
A_i$, $\forall i\in I$.

First suppose that $N\subset A^{\sharp}$. Assume that the total
intersection of $\mathcal{F}$ is empty. Then $R/A$ is non-complete in
its ideal topology. By Proposition~\ref{P:breath}
there exists $x\in\widehat{R}\setminus R$ such that
$\mathrm{B}(x)=A$. Let $b\in A^{\sharp}\setminus N$. There exists
$a\in (A:b)\setminus A$. Since $\mathrm{B}(x)=A$ we have $x=r+ay$ for some
$r\in R$ and $y\in\widehat{R}$. By Lemma~\ref{L:breath}
$\mathrm{B}(y)=(A:a)$. Since $b\in\mathrm{B}(y)$ we have
$N\subset\mathrm{B}(y)$. By Proposition~\ref{P:breath}
$R/\mathrm{B}(y)$ is non-complete in its ideal topology. This
contradicts that $R/N$ is maximal. So the total intersection of
$\mathcal{F}$ is non-empty in this case.

Now we assume that $N=A^{\sharp}$. Then $A$ is an ideal of $R_N$. By
\cite[Lemma 29]{Couch03} either $A=Na$ for some $a\in N$ or
$A=\cap_{a\notin A}aR_N$.

First we assume that $A=Na$. We may suppose that $A_i\subseteq aR_N$,
$\forall i\in I$. Since $\mathcal{F}$ has the finite intersection
property, $x_i+aR_N=x_j+aR_N$, $\forall i,j\in I$. Let $y\in x_i+aR_N$
for each $i\in I$. Then $(x_i-y+A_i)_{i\in I}$ is a family of cosets
of $aR_N$ which has the finite intersection property. But $aR_N/aN$ is
a uniserial module over $R/N$. Then $aR_N/aN$ is linearly compact
since $R/N$ is maximal. Thus $\cap_{i\in
  I}(x_i-y+A_i)\ne\emptyset$. Hence the total intersection of
$\mathcal{F}$ is non-empty.

Now suppose that $A=\cap_{a\notin A}aR_N$. By Proposition~\ref{P:countable} and
\cite[Lemma 30]{Couch03} there exists a countable family
$(a_n)_{n\in\mathbb{N}}$ of elements of $N\setminus A$ such that
$A=\cap_{n\in\mathbb{N}}a_nR_N$ and $a_n\notin a_{n+1}R_N$, $\forall
n\in\mathbb{N}$. By induction on $n$ we get a subfamily
$(A_{i_n})_{n\in\mathbb{N}}$ of the family $(A_i)_{i\in I}$ such that
$A_{i_n}\subset a_nR_N$ in the following way: we choose $i_0\in I$
such that $A_{i_0}\subset a_0R_N$ and, $\forall n\in\mathbb{N}$, we
pick $i_{n+1}$ such that $A_{i_{n+1}}\subset A_{i_n}\cap
a_{n+1}R_N$. Then the family $(x_{i_n}+a_nR_N)_{n\in\mathbb{N}}$ has 
the finite intersection property. Since $R_N$ is maximal there exists
$x\in x_{i_n}+a_nR_N$, $\forall n\in\mathbb{N}$. But the equality
$A=\cap_{a\notin A}aR_N$ implies that, $\forall n\in\mathbb{N}$, there
exists an integer $m>n$ such that $a_mR_N\subseteq A_{i_n}$. Since
$x-x_{i_m}\in a_mR_N$ and $x_{i_m}-x_{i_n}\in A_{i_n}$ we get that
$x\in x_{i_n}+A_{i_n}$, $\forall n\in\mathbb{N}$. Hence $\mathcal{F}$
has a non-empty total intersection. The proof is now complete.
\qed
\begin{theorem} \label{T:malcev}  Then $R$ is maximal if
  and only if, for each polyserial $R$-module $M$, $\widehat{M}$ is
  direct sum of $\mathrm{Mr}\ M$ indecomposable pure-injective modules.  
\end{theorem}
\textbf{Proof.} If $R$ is maximal, then each polyserial module $M$ is a
direct sum of $\mathrm{Mr}\ M$ pure-injective uniserial modules by
\cite[Proposition XII.2.4]{FuSa01} (even if $R$ is not a domain, this proposition holds, with the same proof). 

If $R$ is not maximal then $R/N$ or $R_N$ is not maximal by Theorem~\ref{T:maximal}.

Assume that $R'=R/N$ is not maximal. Then
$E=\widehat{R'}/R'$ is a nonzero torsion-free $R'$-module. Let
$x\in\widehat{R'}\setminus R'$, $\bar{x}$ be its image in $E$ and $U$
the submodule of $E$ such that $U/R'\bar{x}$ is the torsion submodule
of $E/R'\bar{x}$. Then $U$ is a pure submodule of $E$, a rank one
torsion-free module and a uniserial module. Let $M$ be the inverse
image of $U$ by the natural map
$\widehat{R'}\rightarrow E$. Then $M$ is a pure
    submodule of $\widehat{R'}$ and a non-uniserial polyserial
    module with the two following (standard) uniserial
    factors: $R'$ and $U$. We have $\mathrm{Mr}\ M=2$. Let $W$ be a submodule of $\widehat{R'}$ such
    that $M\cap W=0$ and $M\rightarrow \widehat{R'}/W$ is a pure
    monomorphism. Thus $R'\cap W=0$ and $R'\rightarrow \widehat{R'}/W$ is
    a pure monomorphism too. Since $R'$ is pure-essential in $\widehat{R'}$ it follows that $W=0$. We conclude that $M$ is pure-essential in $\widehat{R'}$, so that
    $\widehat{M}=\widehat{R'}\subset\widehat{R'}\oplus\widehat{U}$. (Let
    us observe that $M$ and $U$ are not finitely generated by
    Theorem~\ref{T:main}.) 

Suppose that $R'=R_N$ is not maximal. After replacing $R'$ with
$R'/rR'$, where $r$ is a non-unit of $R'$, we may assume that $R'$ is
coherent and self fp-injective by \cite[Theorem 11]{Couch03}. Then
$E=\widehat{R'}/R'$ is a nonzero fp-injective $R'$-module. By
\cite[Lemma 6]{Cou05} $E$ contains a pure uniserial submodule $U$. We
define $M$ as above. Then $\mathrm{Mr}\ M=2$ and $M$ is an essential
submodule of $\widehat{R'}$. So $\widehat{M}=\widehat{R'}$.
\qed

\end{document}